\documentclass[11pt]{article}
\usepackage{latexsym}
\newcommand{\eproof}{\mbox{\ }\hfill $\Box$ \par \vskip 10pt}

\newtheorem{Theorem}{Theorem}[section]
\newtheorem{lemma}[Theorem]{Lemma}
\newtheorem{prop}[Theorem]{Proposition}

\begin{document}

\title{High frequency dispersive estimates  in dimension two}

\author{{\sc Simon Moulin}}

\date{}

\maketitle

\abstract{We prove dispersive estimates at high frequency in
dimension two for both the wave and the Schr\"odinger groups
for a very large class of real-valued potentials. 

\setcounter{section}{0}
\section{Introduction and statement of results}

The purpose of this note is to prove dispersive estimates at high frequency for the wave group $e^{it\sqrt{G}}$ and the
Schr\"odinger group $e^{itG}$, where $G$ denotes the self-adjoint realization of the
operator $-\Delta+V$ on $L^2({\bf R}^2)$ and $V$ is a real-valued potential which decays at infinity in a way that $G$
has no real resonances nor eigenvalues in an interval $[a_0,+\infty)$, $a_0>0$. In fact, we are looking for as large as
possible class of potentials for which we have dispersive estimates similar to those we do for the free operator
$G_0$. Hereafter $G_0$ denotes the self-adjoint realization of the operator $-\Delta$ on $L^2({\bf R}^2)$.
It turns out that in dimension two one can get such dispersive estimates at high frequency for potentials satisfying
$$\sup_{y\in {\bf R}^2}\int_{{\bf R}^2}\frac{|V(x)|dx}{|x-y|^{1/2}}\le C<+\infty.\eqno{(1.1)}$$
Clearly, (1.1) is fulfilled for potentials $V\in L^\infty({\bf R}^2)$ satisfying 
$$|V(x)|\le C\langle x\rangle^{-\delta},\quad\forall x\in {\bf
R}^2,\eqno{(1.2)}$$ with constants $C>0$, $\delta>3/2$. Given any $a>0$, set
$\chi_a(\sigma)=\chi_1(\sigma/a)$, where $\chi_1\in C^\infty({\bf
R})$, $\chi_1(\sigma)=0$ for $\sigma\le 1$, $\chi_1(\sigma)=1$ for
$\sigma\ge 2$. Our first result is the following

\begin{Theorem} Let $V$ satisfy (1.1). Then, there exists a constant $a_0>0$ so that
for every $a\ge a_0$, $0<\epsilon\ll 1$, $2\le p<+\infty$, we have the estimates
$$\left\|e^{it\sqrt{G}}G^{-3/4-\epsilon}\chi_a(G)\right\|_{L^1\to L^\infty}\le
C_\epsilon|t|^{-1/2},\quad t\neq 0,\eqno{(1.3)}$$ 
$$\left\|e^{it\sqrt{G}}G^{-3\alpha /4}\chi_a(G)
\right\|_{L^{p'}\to L^p}\le
C|t|^{-\alpha/2},\quad t\neq 0,\eqno{(1.4)}$$ where $1/p+1/p'=1$,
$\alpha=1-2/p$. 
\end{Theorem}

The estimate (1.3) is proved in \cite{kn:CCV} under tha assumption (1.2). Moreover, if in addition one supposes that $G$ has no
strictly positive resonances, it is shown in \cite{kn:CCV} that (1.3) holds for any $a>0$ still under (1.2). In dimension three an
analogue of (1.3) is proved in \cite{kn:CCV}, \cite{kn:GV} for potentials satisfying (1.2) with $\delta>2$, and extended in \cite{kn:DP}
to a large subset of potentials satisfying 
$$\sup_{y\in {\bf R}^3}\int_{{\bf R}^3}\frac{|V(x)|dx}{|x-y|}\le C<+\infty.\eqno{(1.5)}$$
In dimensions $n\ge 4$ there are very few results. In \cite{kn:B}, an analogue of (1.3) is proved for potentials belonging to the 
Schwartz class, while in \cite{kn:V1} dispersive estimates with a loss of $(n-3)/2$ derivatives are obtained for potentials satisfying (1.2) with $\delta>(n+1)/2$. Recently, in \cite{kn:M} dispersive estimates at low frequency have been proved in dimensions $n\ge 4$
for a very large class of potentials, provided zero is neighter an eigenvalue nor a resonance. 

Our second result is the following

\begin{Theorem} Let $V$ satisfy (1.1). Then, there exists a constant $a_0>0$ so that
for every $a\ge a_0$, we have the estimate
$$\left\|e^{itG}\chi_a(G)\right\|_{L^1\to L^\infty}\le
C|t|^{-1},\quad t\neq 0.\eqno{(1.6)}$$ 
\end{Theorem}

Note that the estimate (1.6) (for any $a>0$) is proved in \cite{kn:S} for potentials satisfying (1.2) with $\delta>2$.
In dimension three an analogue of (1.6) (for any $a>0$) is proved in \cite{kn:RS} for potentials satisfying (1.5) with $C>0$ small enough,
and in \cite{kn:G} for potentials $V\in L^{3/2-\epsilon}\cap L^{3/2+\epsilon}$, $0<\epsilon\ll 1$, not necessarily small.
In dimensions $n\ge 4$, an analogue of (1.6) (for any $a>0$) is proved in \cite{kn:JSS} for potentials satisfying (1.2) with $\delta>n$
as well as the condition $\widehat V\in L^1$. This result has been recently extended in \cite{kn:MV} to potentials satisfying (1.2) with $\delta>n-1$ and $\widehat V\in L^1$.
Note also the work \cite{kn:V2}, where an analogue of (1.6) (for any $a>0$) with a loss
of $(n-3)/2$ derivatives is obtained for potentials satisfying (1.2) with $\delta>(n+2)/2$.
In \cite{kn:MV} dispersive estimates at low frequency have been also proved in dimensions $n\ge 4$
for a very large class of potentials, provided zero is neighter an eigenvalue nor a resonance. 

To prove (1.3) we use the same idea we have already used in \cite{kn:M} to prove low frequency dispersive estimates in dimensions $n\ge 4$.
The key point is the following estimate which holds in all dimensions $n\ge 2$:
$$h\int_{-\infty}^\infty\left\|Ve^{it\sqrt{G_0}}\psi(h^2G_0)f\right\|_{L^1}dt\le \gamma_n C_n(V)h^{-(n-3)/2}\left\|f\right\|_{L^1},\quad h>0,\eqno{(1.7)}$$ 
where $\psi\in C_0^\infty((0,+\infty))$, $\gamma_n>0$ is a constant independent of $V$, $h$ and $f$, and
$$C_n(V):=\sup_{y\in {\bf R}^n}\int_{{\bf R}^n}\frac{|V(x)|dx}{|x-y|^{(n-1)/2}}<+\infty.\eqno{(1.8)}$$
Our approach is based on the observation that if
$$C_n(V)h^{-(n-3)/2}\ll 1,\eqno{(1.9)}$$
then (1.7) implies (under reasonable assumptions on the potential) a similar estimate for the perturbed wave group, namely
$$h\int_{-\infty}^\infty\left\|Ve^{it\sqrt{G}}\psi(h^2G)f\right\|_{L^1}dt\le \widetilde C_n(V)h^{-(n-3)/2}\left\|f\right\|_{L^1}. \eqno{(1.10)}$$ 
When $n=3$, (1.9) is fulfilled for small potentials and all $h$, when $n\ge 4$, (1.9) is fulfilled for large $h$ (i.e. at low frequency)
without extra restrictions on the potential, while for $n=2$, (1.9) is fulfilled for small $h$ (i.e. at high frequency) again
without restrictions on the potential others than (1.1). Note that (1.10) may hold without (1.9). Indeed, when $n=3$, (1.10)
is proved in \cite{kn:G} for potentials $V\in L^{3/2-\epsilon}\cap L^{3/2+\epsilon}$ and all $h>0$, and then used to prove the three
dimensional analogue of (1.6). In the present paper we adapt this approach to the case of dimension two, and show that (1.6) follows
from (1.10) for potentials satisfying (1.1) only, provided the parameter $a$ is taken large enough (see Section 3).

\section{Proof of Theorem 1.1}

Let $\psi\in C_0^\infty((0,+\infty))$ and set
$$\Phi(t;h)=e^{it\sqrt{G}}\psi(h^2G)-e^{it\sqrt{G_0}}\psi(h^2G_0).$$
We will first show that (1.3) and (1.4) follow
 from the following

\begin{prop} Let $V$ satisfy (1.1). Then, there exist
positive constants $C$ and $h_0$ so that for $0<h\le h_0$ we
have
$$\left\|\Phi(t;h)\right\|_{L^1\to L^\infty}\le
Ch^{-1}|t|^{-1/2},\quad t\neq 0.\eqno{(2.1)}$$
\end{prop}

Writing
$$\sigma^{-3/4-\epsilon}\chi_a(\sigma)=\int_0^{a^{-1}}\psi(\sigma\theta)
\theta^{-1/4+\epsilon}d\theta, $$
where $\psi(\sigma)=\sigma^{1/4-\epsilon}\chi'_1(\sigma)\in
C_0^\infty((0,+\infty))$, and using (2.1) we get 
$$\left\|e^{it\sqrt{G}}G^{-3/4-\epsilon}\chi_a(G)-e^{it\sqrt{G_0}}
G_0^{-3/4-\epsilon}\chi_a(G_0)\right\|_{L^1\to L^\infty}$$ $$\le
\int_0^{a^{-1}}\left\|\Phi(t;\sqrt{\theta})\right\|_{L^1\to
L^\infty}\theta^{-1/4+\epsilon}d\theta\le 
C|t|^{-1/2}\int_0^{a^{-1}}\theta^{-3/4+\epsilon}d\theta
\le C|t|^{-1/2},\eqno{(2.2)}$$ provided $a$ is taken large 
enough. Clearly, (1.3) follows from (2.2) and the fact that it holds
for $G_0$. To prove (1.4), observe that an interpolation between (2.1) and the trivial bound
$$\left\|\Phi(t;h)\right\|_{L^2\to L^2}\le C$$
yields
$$\left\|\Phi(t;h)\right\|_{L^{p'}\to L^p}\le
Ch^{-\alpha}|t|^{-\alpha/2},\quad t\neq 0,\eqno{(2.3)}$$
for every $2\le p\le +\infty$, $p'$ and $\alpha$ being as in Theorem 1.1.  Now we write
$$\sigma^{-3\alpha/4}\chi_a(\sigma)=\int_0^{a^{-1}}\psi(\sigma\theta)
\theta^{-1+3\alpha/4}d\theta, $$
and use (2.3) to obtain (for $0<\alpha\le 1$)
$$\left\|e^{it\sqrt{G}}G^{-3\alpha/4}\chi_a(G)-e^{it\sqrt{G_0}}
G_0^{-3\alpha/4}\chi_a(G_0)\right\|_{L^{p'}\to L^p}$$ $$\le
\int_0^{a^{-1}}\left\|\Phi(t;\sqrt{\theta})\right\|_{L^{p'}\to
L^p}\theta^{-1+3\alpha/4}d\theta\le 
C|t|^{-\alpha/2}\int_0^{a^{-1}}\theta^{-1+\alpha/4}d\theta
\le C|t|^{-\alpha/2},\eqno{(2.4)}$$ provided $a$ is taken large 
enough. Now, (1.4) follows from (2.4) and the fact that it holds
for $G_0$.\\

{\it Proof of Proposition 2.1.} We will first prove the following

\begin{lemma} Let $V$ satisfy (1.1). Then, there exist
positive constants $C$ and $h_0$ so that for $0<h\le h_0$ we
have
$$\left\|\psi(h^2G)-\psi(h^2G_0)\right\|_{L^1\to L^1}\le Ch^{1/2}.\eqno{(2.5)}$$
\end{lemma}

{\it Proof.} We will make use of the formula
$$\psi(h^2G)=\frac{2}{\pi}\int_{\bf C}\frac{\partial\widetilde\varphi}{\partial\bar z}(z)(h^2G-z^2)^{-1}zL(dz),\eqno{(2.6)}$$
where $L(dz)$ denotes the Lebesgue measure on ${\bf C}$, $\widetilde\varphi\in C_0^\infty({\bf C})$ is an almost analytic continuation of
$\varphi(\lambda)=\psi(\lambda^2)$ supported in a small complex neighbourhood of supp$\,\varphi$ and satisfying
$$\left|\frac{\partial\widetilde\varphi}{\partial\bar z}(z)\right|\le C_N|{\rm Im}\,z|^N,\quad\forall N\ge 1.$$
For $\pm {\rm Im}\,\lambda\ge 0$, ${\rm Re}\,\lambda> 0$, set
$$R_0^\pm(\lambda)=(G_0-\lambda^2)^{-1},\quad R^\pm(\lambda)=(G-\lambda^2)^{-1}.$$
We have the identity
$$R^\pm(\lambda)\left(1+VR_0^\pm(\lambda)\right)=R_0^\pm(\lambda).\eqno{(2.7)}$$
It is well known that the kernels of the operators $R_0^\pm(\lambda)$ are given in terms of the zero order Hankel functions
by the formula
$$[R_0^\pm(\lambda)](x,y)=\pm i4^{-1}H_0^\pm(\lambda|x-y|).$$
Moreover, the functions $H_0^\pm$ satisfy the bound
$$\left|H_0^\pm(\lambda)\right|\le C|\lambda|^{-1/2}e^{-|{\rm Im}\,\lambda|},\quad|\lambda|\ge 1,\,\pm {\rm Im}\,\lambda\ge 0,\eqno{(2.8)}$$
while near $\lambda=0$ they are of the form
$$H_0^\pm(\lambda)=a_{0,1}^\pm(\lambda)+a_{0,2}^\pm(\lambda)\log\lambda,\eqno{(2.9)}$$
where $a_{0,j}^\pm$ are analytic functions. In particular, we have
$$\left|H_0^\pm(\lambda)\right|\le C|\lambda|^{-1/2},\quad{\rm Re}\,\lambda>0,\,\pm {\rm Im}\,\lambda\ge 0.\eqno{(2.10)}$$
Using these bounds we will prove the following

\begin{lemma} Let $V$ satisfy (1.1). Then, there exist constants $C>0$ and $0<h_0\le 1$ so that for $z\in {\bf C}_\varphi^\pm:=\{z\in{\rm supp}\,\widetilde\varphi,\,\pm{\rm Im}\,z\ge 0\}$, we have the estimates
$$\left\|VR_0^\pm(z/h)\right\|_{L^1\to L^1}\le Ch^{1/2},\quad 0<h\le 1,\eqno{(2.11)}$$
$$\left\|VR^\pm(z/h)\right\|_{L^1\to L^1}\le Ch^{1/2},\quad 0<h\le h_0,\eqno{(2.12)}$$
$$\left\|R_0^\pm(z/h)\right\|_{L^1\to L^1}\le Ch^2|{\rm Im}\,z|^{-2},\quad 0<h\le 1,\,{\rm Im}\,z\neq 0,\eqno{(2.13)}$$
$$\left\|R^\pm(z/h)\right\|_{L^1\to L^1}\le Ch^2|{\rm Im}\,z|^{-2},\quad 0<h\le h_0,\,{\rm Im}\,z\neq 0.\eqno{(2.14)}$$
\end{lemma}

{\it Proof.} By (1.1) and (2.10), the norm in the LHS of (2.11) is upper bounded by
$$\sup_{y\in{\bf R}^2}\int_{{\bf R}^2}|V(x)||H_0^\pm(z|x-y|/h)|dx\le  Ch^{1/2}\int_{{\bf R}^2}\frac{|V(x)dx}{|x-y|^{1/2}}\le C'h^{1/2}.$$
Similarly, the norm in the LHS of (2.13) is upper bounded by
$$\sup_{y\in{\bf R}^2}\int_{{\bf R}^2}|H_0^\pm(z|x-y|/h)|dx=h^2\sup_{y\in{\bf R}^2}\int_{{\bf R}^2}|H_0^\pm(z|x-y|)|dx$$ 
$$\le Ch^2|{\rm Im}\,z|^{-2}\int_{{\bf R}^2}\langle x-y\rangle^{-3/2}|x-y|^{-1}dx=C'h^2|{\rm Im}\,z|^{-2}\int_0^\infty
\langle\sigma\rangle^{-3/2}d\sigma.$$
To prove (2.12) and (2.14) we will make use of the identity (2.7). It follows from (2.11) that there exists a constant 
$0<h_0\le 1$ so that for $0<h\le h_0$ the operator $1+VR_0^\pm(z/h)$ is invertible on $L^1$ with an inverse satisfying
$$\left\|\left(1+VR_0^\pm(z/h)\right)^{-1}\right\|_{L^1\to L^1}\le C,\quad z\in {\bf C}_\varphi^\pm,\eqno{(2.15)}$$
with a constant $C>0$ independent of $z$ and $h$. Clearly, (2.12) follows from (2.11) and (2.15), while (2.14) follows from (2.13)
and (2.15).
\eproof

To prove (2.5) we rewrite the identity (2.7) in the form
$$R^\pm(z/h)-R_0^\pm(z/h)=R_0^\pm(z/h)VR_0^\pm(z/h)\left(1+VR_0^\pm(z/h)\right)^{-1},$$
and hence, using Lemma 2.3 and (2.15), we get
$$\left\|h^{-2}R^\pm(z/h)-h^{-2}R_0^\pm(z/h)\right\|_{L^1\to L^1}\le Ch^{1/2}|{\rm Im}\,z|^{-2},\quad 0<h\le h_0,\,
z\in {\bf C}_\varphi^\pm,\,{\rm Im}\,z\neq 0.\eqno{(2.16)}$$
It is easy now to see that (2.5) follows from (2.6) and (2.16).
\eproof

We will now derive (2.1) from the
following

\begin{prop} Let $V$ satisfy (1.1). Then, there exist
positive constants $C$ and $h_0$ so that we have, for $0\le
s\le 1/2$, $f,g\in L^1$,
$$\left\|e^{it\sqrt{G_0}}\psi(h^2G_0)f\right\|_{L^\infty}\le
Ch^{-3/2}|t|^{-1/2}\|f\|_{L^1},\quad h>0,\,t\neq
0,\eqno{(2.17)}$$
$$\int_{{\bf R}^2}\int_{{\bf R}^2}\int_{-\infty}^\infty
|t|^s|x-y|^{-s}
\left|Ve^{it\sqrt{G_0}}\psi(h^2G_0)f(x)
\right|\left|g(y)\right|dt dxdy$$ $$\le Ch^{-1/2}
\|f\|_{L^1}\|g\|_{L^1},\quad h>0,\eqno{(2.18)}$$
$$\int_{{\bf R}^2}\int_{{\bf R}^2}\int_{-\infty}^\infty
|t|^s\left\langle|x-y|/h\right\rangle^{-s}
\left|Ve^{it\sqrt{G}}\psi(h^2G)f(x)
\right|\left|g(y)\right|dt dxdy$$ $$\le Ch^{s-1/2}
\|f\|_{L^1}\|g\|_{L^1},\quad 0<h\le h_0.\eqno{(2.19)}$$
\end{prop}

As in \cite{kn:V2}, using Duhamel's formula
$$e^{it\sqrt{G}}=e^{it\sqrt{G_0}}+i\frac{\sin\left(t\sqrt{G_0}\right)}
{\sqrt{G_0}}
\left(\sqrt{G}-\sqrt{G_0}\right) -\int_0^t
\frac{\sin\left((t-\tau)\sqrt{G_0}\right)}{\sqrt{G_0}}Ve^{i\tau
\sqrt{G}}d\tau$$ we get the identity
$$\Phi(t;h)=\sum_{j=1}^2\Phi_j(t;h),\eqno{(2.20)}$$
where
$$\Phi_1(t;h)=\left(\psi_1(h^2G)-\psi_1(h^2G_0)\right)e^{it\sqrt{G}}
\psi(h^2G)
$$ $$+\psi_1(h^2G_0)e^{it\sqrt{G_0}}\left(\psi(h^2G)-\psi(h^2G_0)\right)$$
 $$-i\psi_1(h^2G_0)\sin\left(t\sqrt{G_0}\right)\left(\psi(h^2G)
-\psi(h^2G_0)\right)
  $$ $$+i\widetilde\psi_1(h^2G_0)\sin\left(t\sqrt{G_0}\right)\left(
 \widetilde\psi(h^2G)-\widetilde\psi(h^2G_0)\right),$$
$$\Phi_2(t;h)=-h\int_0^t\widetilde\psi_1(h^2G_0)\sin\left((t-\tau)
\sqrt{G_0}
\right)Ve^{i\tau\sqrt{ G}}\psi(h^2G)d\tau,$$ where $\psi_1\in
C_0^\infty((0,+\infty))$, $\psi_1=1$ on supp$\,\psi$,
$\widetilde\psi(\sigma)=\sigma^{1/2}\psi(\sigma)$,
$\widetilde\psi_1(\sigma)=\sigma^{-1/2}\psi_1(\sigma)$. By
Proposition 2.4 and (2.5), we have
$$\left\|\Phi_1(t;h)f\right\|_{L^\infty}\le
Ch^{-1}|t|^{-1/2}\|f\|_{L^1}+Ch^{1/2}
\left\|\Phi(t;h)f\right\|_{L^\infty},\eqno{(2.21)}$$
$$t^{1/2}\left|\langle\Phi_2(t;h)f,g\rangle\right|$$ $$\le h
\int_0^{t/2}(t-\tau)^{1/2} \left\|\sin\left((t-\tau)\sqrt{G_0}
\right)\widetilde\psi_1(h^2G_0)g\right\|_{
L^\infty}\left\|Ve^{i\tau \sqrt{G}}\psi(h^2G)f\right\|_{L^1}d\tau$$
 $$+h\int_{t/2}^t\left\|V\sin\left((t-\tau)\sqrt{G_0}
\right)\widetilde\psi_1(h^2G_0)g\right\|_{ L^1} \tau^{1/2}
\left\|e^{i\tau \sqrt{G}}\psi(h^2G)f\right\|_{L^\infty}d\tau$$
 $$\le Ch^{-1/2}\|g\|_{ L^1}
\int_{-\infty}^{\infty}\left\|Ve^{i\tau\sqrt{
G}}\psi(h^2G)f\right\|_{L^1}d\tau$$
 $$+h\sup_{t/2\le\tau\le t}\tau^{1/2} \left\|e^{i\tau\sqrt{
G}}\psi(h^2G)f\right\|_{L^\infty}
\int_{-\infty}^{\infty}\left\|V\sin\left((t-\tau)\sqrt{G_0}
\right)\widetilde\psi_1(h^2G_0)g\right\|_{L^1}d\tau$$
 $$\le Ch^{-1}\|g\|_{L^1}\|f\|_{L^1}
 +Ch^{1/2}\|g\|_{L^1}\sup_{t/2\le\tau\le t}\tau^{1/2} 
\left\|e^{i\tau\sqrt{G}}\psi(h^2G)f\right\|_{L^\infty},$$ for $t>0$, which clearly implies
 $$t^{1/2}\left\|\Phi_2(t;h)f\right\|_{L^\infty}\le
Ch^{-1}\|f\|_{L^1}+Ch^{1/2}\sup_{t/2\le\tau\le
t}\tau^{1/2} \left\|e^{i\tau\sqrt{
G}}\psi(h^2G)f\right\|_{L^\infty}.\eqno{(2.22)}$$ By (2.20)-(2.22),
we conclude
$$t^{1/2}\left\|\Phi(t;h)f\right\|_{L^\infty}\le
Ch^{-1}\|f\|_{L^1}+Ch^{1/2}
t^{1/2}\left\|\Phi(t;h)f\right\|_{L^\infty}$$ $$
+Ch^{1/2}\sup_{t/2\le\tau\le t}\tau^{1/2}
\left\|\Phi(\tau;h)f\right\|_{L^\infty}.\eqno{(2.23)}$$ Taking $h$
small enough we can absorb the second and the third terms in the RHS
of (2.23), thus obtaining (2.1). Clearly, the case of $t<0$ can be
treated in the same way.\eproof

{\it Proof of Proposition 2.3.} The kernel of the operator $e^{it\sqrt{ G_0}}\psi(h^2G_0)$ is of
the form $K_h(|x-y|,t)$, where
$$K_h(\sigma,t)=(2\pi)^{-1}\int_0^\infty
e^{it\lambda}J_0(\sigma\lambda)\psi(h^2\lambda^2)\lambda
d\lambda= h^{-2}K_1(\sigma h^{-1},th^{-1}),\eqno{(2.24)}$$ where
$J_0(z)=\left(H_0^+(z)+H_0^-(z)\right)/2$ is the Bessel
function of order zero. It is shown in \cite{kn:V2}
(Section 2) that $K_h$ satisfies the estimates (for all $\sigma,
h>0$, $t\neq 0$)
$$\left|K_1(\sigma,t)\right|\le
C|t|^{-s}\langle\sigma\rangle^{s-1/2},\quad\forall s\ge
0,\eqno{(2.25)}$$
$$\left|K_h(\sigma,t)\right|\le
Ch^{-3/2}|t|^{-s}\sigma^{s-1/2},\quad 0\le s\le 1/2.\eqno{(2.26)}$$ Clearly, (2.17) follows from (2.26) with
$s=1/2$. It is easy also to see that (2.18) follows from (1.1) and the following

\begin{lemma} For all $0\le s\le 1/2$, $\sigma,h>0$, we have
 $$\int_{-\infty}^\infty\langle t/h\rangle^s\left|K_h(\sigma,t)\right|dt\le
Ch^{-1}\langle\sigma/h\rangle^{s-1/2},\eqno{(2.27)}$$
$$\int_{-\infty}^\infty|t|^s\left|K_h(\sigma,t)\right|dt\le
Ch^{-1/2}\sigma^{s-1/2}.\eqno{(2.28)}$$
\end{lemma}

{\it Proof.} Clearly, (2.28) follows from (2.27). It is also clear from (2.24) that it suffices to prove (2.27)
with $h=1$. When $0<\sigma\le 1$, this follows from (2.25). Let now $\sigma\ge 1$. We decompose $K_1$ as $K_1^++K_1^-$,
where $K_1^\pm$ is defined by replacing in (2.24) the function $J_0$ by $H_0^\pm/2$. Recall that $H_0^\pm(z)=e^{\pm iz}b_0^\pm(z)$,
where $b_0^\pm(z)$ is a symbol of order $-1/2$ for $z\ge 1$. Using this fact and integrating by parts $m$ times, we get
$$|K_1^\pm(\sigma,t)|\le C_m\sigma^{-1/2}|t\pm\sigma|^{-m}.\eqno{(2.29)}$$
By (2.29), we obtain
$$\int_{-\infty}^\infty\langle t\rangle^s\left|K_1^\pm(\sigma,t)\right|dt\le 2\sigma^s\int_{-\infty}^\infty\left|K_1^\pm(\sigma,t)\right|dt
+\int_{-\infty}^\infty|t\pm\sigma|^s\left|K_1^\pm(\sigma,t)\right|dt$$
 $$\le C'_m\sigma^{s-1/2}\int_{-\infty}^\infty |t\pm\sigma|^{-m}dt
+C_m\sigma^{-1/2}\int_{-\infty}^\infty|t\pm\sigma|^{s-m}dt\le C\sigma^{s-1/2},$$
which clearly implies (2.27) in this case.
\eproof

To prove (2.19) we will use the formula
$$e^{it\sqrt{ G}}\psi(h^2G)=(i\pi h)^{-1}\int_0^\infty e^{it\lambda}\varphi_h(\lambda)\left(R^+(\lambda)-R^-(\lambda)\right)d\lambda,\eqno{(2.30)}$$
where $\varphi_h(\lambda)=\varphi_1(h\lambda)$, $\varphi_1(\lambda)=\lambda\psi(\lambda^2)$.
Combining (2.30) together with (2.7), we get
$$Ve^{it\sqrt{ G}}\psi(h^2G)=(i\pi h)^{-1}\sum_\pm\pm\int_{-\infty}^\infty VP_h^\pm(t-\tau)U_h^\pm(\tau)d\tau,\eqno{(2.31)}$$
where
$$P_h^\pm(t)=\int_0^\infty e^{it\lambda}\widetilde\varphi_h(\lambda)R^\pm_0(\lambda)d\lambda,$$
$$U_h^\pm(t)=\int_0^\infty e^{it\lambda}\varphi_h(\lambda)\left(1+VR^\pm_0(\lambda)\right)^{-1}d\lambda,$$
where $\widetilde\varphi_h(\lambda)=\widetilde\varphi_1(h\lambda)$, $\widetilde\varphi_1\in C_0^\infty((0,+\infty))$ is such that
$\widetilde\varphi_1=1$ on supp$\,\varphi_1$. The kernel of the operator $P_h^\pm(t)$ is of the form $A_h^\pm(|x-y|,t)$, where
$$A_h^\pm(\sigma,t)=\pm i4^{-1}\int_0^\infty e^{it\lambda}\widetilde\varphi_h(\lambda)H^\pm_0(\sigma\lambda)d\lambda=
h^{-1}A_1^\pm(\sigma/h,t/h).\eqno{(2.32)}$$
In the same way as in the proof of Lemma 2.5 one can see that the function $A_h^\pm$ satisfies the estimate
$$\int_{-\infty}^\infty|t|^s\left|A_h^\pm(\sigma,t)\right|dt\le
Ch^{1/2}\sigma^{s-1/2}(1+h^{\epsilon_s}\sigma^{-\epsilon_s}),\quad 0\le s\le 1/2,\quad 0<h\le 1,\eqno{(2.33)}$$
where $\epsilon_s=0$ if $0\le s<1/2$, $\epsilon_s=\epsilon$ if $s=1/2$.

Clearly, it suffices to prove (2.19) with $s=0$ and $s=1/2$. For these values of $s$, 
using (1.1), (2.31) and (2.33), we obtain
$$\int_{{\bf R}^2}\int_{{\bf R}^2}\int_{-\infty}^\infty|t|^s\left\langle|x-y|/h\right\rangle^{-s}
\left|Ve^{it\sqrt{G}}\psi(h^2G)f(x)\right|\left|g(y)\right|dt dxdy$$ $$\le Ch^{-1}\sum_\pm\int_{{\bf R}^2}\int_{{\bf R}^2}\int_{-\infty}^\infty\int_{-\infty}^\infty  \left\langle|x-y|/h\right\rangle^{-s}\left(|t-\tau|^s+|\tau|^s\right)$$ $$\times
\left|VP_h^\pm(t-\tau)U_h^\pm(\tau)f(x)\right|\left|g(y)\right|d\tau dt dxdy$$
$$\le Ch^{-1}\sum_\pm\int_{{\bf R}^2}\int_{{\bf R}^2}\int_{{\bf R}^2}\int_{-\infty}^\infty\int_{-\infty}^\infty  |V(x)|
\left\langle|x-y|/h\right\rangle^{-s}\left(|t-\tau|^s+|\tau|^s\right)$$ $$\times
|A_h^\pm(|x-x'|,t-\tau)|\left|U_h^\pm(\tau)f(x')\right|\left|g(y)\right|d\tau dt dx'dxdy$$
$$\le Ch^{-1}\sum_\pm\int_{{\bf R}^2}\int_{{\bf R}^2}\int_{{\bf R}^2} |V(x)|\left\langle|x-y|/h\right\rangle^{-s}\left|g(y)\right|$$ $$\times\left(\int_{-\infty}^\infty|\tau|^s
|A_h^\pm(|x-x'|,\tau)|d\tau\right)\left(\int_{-\infty}^\infty\left|U_h^\pm(\tau)f(x')\right|d\tau\right) dx'dxdy$$
$$+ Ch^{-1}\sum_\pm\int_{{\bf R}^2}\int_{{\bf R}^2}\int_{{\bf R}^2} |V(x)|\left\langle|x-y|/h\right\rangle^{-s}\left|g(y)\right|$$ $$\times\left(\int_{-\infty}^\infty
|A_h^\pm(|x-x'|,\tau)|d\tau\right)\left(\int_{-\infty}^\infty|\tau|^s\left|U_h^\pm(\tau)f(x')\right|d\tau\right) dx'dxdy$$
 $$\le Ch^{-1/2}\sum_\pm\int_{{\bf R}^2}\int_{{\bf R}^2}\int_{{\bf R}^2} |V(x)|\left\langle|x-y|/h\right\rangle^{-s}|x-x'|^{s-1/2}\left(1+h^{\epsilon_s}|x-x'|^{-\epsilon_s}\right)\left|g(y)\right|$$ $$\times\left(\int_{-\infty}^\infty\left|U_h^\pm(\tau)f(x')\right|d\tau\right) dx'dxdy$$
$$+ Ch^{-1/2}\sum_\pm\int_{{\bf R}^2}\int_{{\bf R}^2}\int_{{\bf R}^2} |V(x)|\left\langle|x-y|/h\right\rangle^{-s}|x-x'|^{-1/2}\left|g(y)\right|$$ $$\times\left(\int_{-\infty}^\infty|\tau|^s\left|U_h^\pm(\tau)f(x')\right|d\tau\right) dx'dxdy:=I_1+I_2.\eqno{(2.34)}$$
To estimate $I_1$ when $s=1/2$, set $q=(2\epsilon)^{-1}$, $1/p+1/q=1$, and observe that in view of (1.1) we have the bound
$$\int_{{\bf R}^2}|V(x)|\left\langle|x-y|/h\right\rangle^{-1/2}|x-x'|^{-\epsilon}dx$$
$$\le\left(\int_{{\bf R}^2}|V(x)|\left\langle|x-y|/h\right\rangle^{-p/2}dx\right)^{1/p}\left(\int_{{\bf R}^2}|V(x)|
|x-x'|^{-1/2}dx\right)^{1/q}$$
$$\le C_1\left(\int_{{\bf R}^2}|V(x)|\left\langle|x-y|/h\right\rangle^{-1/2}dx\right)^{1/p}$$
$$\le C_1h^{1/(2p)}\left(\int_{{\bf R}^2}|V(x)||x-y|^{-1/2}dx\right)^{1/p}\le C_2h^{1/2-\epsilon}.$$
Thus, we obtain
$$I_1\le C'h^{s-1/2}\sum_\pm\int_{{\bf R}^2}\int_{{\bf R}^2}\int_{-\infty}^\infty\left|U_h^\pm(\tau)f(x')\right|\left|g(y)\right| d\tau dx'dy.\eqno{(2.35)}$$
To estimate $I_2$ when $s=1/2$, we use the inequality
$$\left\langle|x-y|/h\right\rangle^{-1/2}|x-x'|^{-1/2}\le \left\langle|x'-y|/h\right\rangle^{-1/2}\left(|x-y|^{-1/2}+|x-x'|^{-1/2}\right).$$
We get
$$I_2\le C''h^{-1/2}\sum_\pm\int_{{\bf R}^2}\int_{{\bf R}^2} \int_{-\infty}^\infty|\tau|^s\left\langle|x'-y|/h\right\rangle^{-s}\left|U_h^\pm(\tau)f(x')\right| \left|g(y)\right|d\tau dx'dy.\eqno{(2.36)}$$
On the other hand, by the identity 
$$\left(1+VR^\pm_0(\lambda)\right)^{-1}=1-VR^\pm_0(\lambda)\left(1+VR^\pm_0(\lambda)\right)^{-1},$$
we obtain
$$U_h^\pm(t)=\widehat\varphi_h(t)-\int_{-\infty}^\infty  VP_h^\pm(t-\tau)U_h^\pm(\tau)d\tau.\eqno{(2.37)}$$
Since
$$\widehat\varphi_h(t)=h^{-1}\widehat\varphi_1(t/h),$$
we have
$$\int_{-\infty}^\infty|t|^s|\widehat\varphi_h(t)|dt\le Ch^s.\eqno{(2.38)}$$
Using (2.37) and (2.38), in the same way as in the proof of (2.34)-(2.36), we obtain with $s=0$ or $s=1/2$,
$$\int_{{\bf R}^2}\int_{{\bf R}^2} \int_{-\infty}^\infty|t|^s\left\langle|x-y|/h\right\rangle^{-s}\left|U_h^\pm(t)f(x)\right|
 \left|g(y)\right|dtdxdy\le Ch^s\|f\|_{L^1}\|g\|_{L^1}$$
 $$+Ch^{s+1/2}\int_{{\bf R}^2}\int_{{\bf R}^2}\int_{-\infty}^\infty\left|U_h^\pm(\tau)f(x')\right|\left|g(y)\right| d\tau dx'dy$$
$$+ Ch^{1/2}\int_{{\bf R}^2}\int_{{\bf R}^2} \int_{-\infty}^\infty|\tau|^s\left\langle|x'-y|/h\right\rangle^{-s}\left|U_h^\pm(\tau)f(x')\right| \left|g(y)\right|d\tau dx'dy.\eqno{(2.39)}$$
Taking $h$ small enough we can absorb the second and the third terms in the RHS of (2.39) and get the estimate
$$\int_{{\bf R}^2}\int_{{\bf R}^2} \int_{-\infty}^\infty|t|^s\left\langle|x-y|/h\right\rangle^{-s}\left|U_h^\pm(t)f(x)\right|
 \left|g(y)\right|dtdxdy\le C'h^s\|f\|_{L^1}\|g\|_{L^1}.\eqno{(2.40)}$$
 Now (2.19) follows from (2.34)-(2.36) and (2.40).
\eproof

\section{Proof of Theorem 1.2}

Set
$$\Psi(t;h)=e^{itG}\psi(h^2G)-e^{itG_0}\psi(h^2G_0).$$
As in the previous section, one
can derive (1.6) from the following

\begin{prop} Let $V$ satisfy (1.1). Then, there exist
positive constants $C$ and $h_0$ so that for $0<h\le h_0$, we
have
$$\left\|\Psi(t;h)\right\|_{L^1\to L^\infty}\le
Ch^{1/2}|t|^{-1},\quad t\neq 0.\eqno{(3.1)}$$
\end{prop}

{\it Proof.} We will derive (3.1) from (2.19). To this end, we will
use the identity
$$e^{it\lambda^2}\varphi(h^2\lambda^2)=\int_{-\infty}^\infty
e^{i\tau\lambda}\zeta_h(t,\tau)d\tau,\eqno{(3.2)}$$ where
$\varphi\in C_0^\infty((0,+\infty))$, $\varphi=1$ on
supp$\,\psi_1$, the functions $\psi$ and $\psi_1$ being as in the previous
section, and
$$\zeta_h(t,\tau)=(2\pi)^{-1}\int_0^\infty
e^{it\lambda^2-i\tau\lambda}\varphi(h^2\lambda^2)d\lambda=h^{-1}
\zeta_1(th^{-2},\tau h^{-1}).\eqno{(3.3)}$$ We deduce from (3.2) the
formula
$$e^{itG}\psi(h^2G)=\int_{-\infty}^\infty
\zeta_h(t,\tau)e^{i\tau\sqrt{G}}\psi(h^2G)d\tau.\eqno{(3.4)}$$
Given any integer $m\ge 0$, integrating by parts $m$ times and using the well known bound
$$\left|\int_{-\infty}^\infty
e^{it\lambda^2-i\tau\lambda}\phi(\lambda)d\lambda\right|\le C|t|^{-1/2},\quad \forall t\neq 0,\,\tau\in{\bf R},$$
where $\phi\in C_0^\infty({\bf R})$, one easily obtains the bound
$$\left|\zeta_1(t,\tau)\right|\le C_m|t|^{-m-1/2}\langle\tau\rangle^m,\quad \forall t\neq 0,\,\tau\in{\bf R}.\eqno{(3.5)}$$
By (3.3) and (3.5),
$$\left|\zeta_h(t,\tau)\right|\le C_m h^{2m}|t|^{-m-1/2}\langle\tau/h\rangle^m,\quad \forall t\neq 0,\,\tau\in{\bf R},\,h>0,\eqno{(3.6)}$$
for every integer $m\ge 0$, and hence for all real $m\ge 0$. 
 By (2.5), (2.20) and (3.4), we get
$$\left|\left\langle\Psi(t;h)f,g\right\rangle\right|\le Ch^{1/2}
\left\|\Psi(t;h)f\right\|_{L^\infty}\|g\|_{L^1}$$
$$+\int_{-\infty}^\infty\left|\zeta_h(t,\tau)\right|
\left|\left\langle e^{i\tau\sqrt{G_0}}\psi(h^2G_0)f,\left(\psi_1(h^2G)-
\psi_1(h^2G_0)\right)g\right\rangle\right|d\tau$$
$$+\int_{-\infty}^\infty\left|\zeta_h(t,\tau)\right|
\left|\left\langle e^{i\tau\sqrt{G_0}}\psi_1(h^2G_0)\left(\psi(h^2G)-
\psi(h^2G_0)\right)f,g\right\rangle\right|d\tau$$
$$+\int_{-\infty}^\infty\left|\zeta_h(t,\tau)\right|
\left|\left\langle \sin\left(\tau\sqrt{G_0}\right)
\psi_1(h^2G_0)\left(\psi(h^2G)-
\psi(h^2G_0)\right)f,g\right\rangle\right|d\tau$$
$$+\int_{-\infty}^\infty\left|\zeta_h(t,\tau)\right|
\left|\left\langle \sin\left(\tau\sqrt{G_0}\right)
\widetilde\psi_1(h^2G_0)\left(\widetilde\psi(h^2G)-
\widetilde\psi(h^2G_0)\right)f,g\right\rangle\right|d\tau$$
$$+h\int_{-\infty}^\infty\int_0^\tau\left|\zeta_h(t,\tau)\right|
\left|\left\langle Ve^{i\tau'\sqrt{G}}\psi(h^2G)f,
 \sin\left((\tau-\tau')\sqrt{G_0}\right)
\widetilde\psi_1(h^2G_0)g\right\rangle\right|d\tau'd\tau.\eqno{(3.7)}$$
Using (3.6) with $m=1/2$ and (2.27) with $s=1/2$ together with (2.5), we obtain that the first integral in the RHS of (3.7) is
upper bounded by
$$Ch|t|^{-1}\int_{{\bf R}^2}\int_{{\bf R}^2}\int_{-\infty}^\infty\langle\tau/h\rangle^{1/2}
\left|K_h(|x-y|,\tau)\right||f(x)|\left|\left(\psi_1(h^2G)-
\psi_1(h^2G_0)\right)g(y)\right|d\tau dx dy$$
 $$\le C|t|^{-1}\int_{{\bf R}^2}\int_{{\bf R}^2}
|f(x)|\left|\left(\psi_1(h^2G)-
\psi_1(h^2G_0)\right)g(y)\right| dx dy\le Ch^{1/2}|t|^{-1}\|f\|_{L^1}\|g\|_{L^1},$$
and similarly for the next three integrals. The last term is upper bounded by
$$Ch^{2}|t|^{-1}\int_{{\bf R}^2}\int_{{\bf R}^2}\int_{-\infty}^\infty
\int_0^\tau\left(|\tau'/h|^{1/2}+\langle(\tau-\tau')/h\rangle^{1/2}\right)\left|\widetilde 
K_h(|x-y|,(\tau-\tau'))\right|$$ $$\times
\left|Ve^{i\tau'\sqrt{G}}\psi(h^2G)f(x)
\right|\left|g(y)\right|d\tau'd\tau dx dy$$
$$\le Ch^{3/2}|t|^{-1}\int_{{\bf R}^2}\int_{{\bf R}^2}\int_{-\infty}^\infty
|\tau|^{1/2}\left|Ve^{i\tau\sqrt{G}}\psi(h^2G)f(x)
\right|\left|g(y)\right|d\tau\int_{-\infty}^\infty\left|\widetilde 
K_h(|x-y|,\tau)\right|d\tau dxdy$$
$$+ Ch^{2}|t|^{-1}\int_{{\bf R}^2}\int_{{\bf R}^2}\int_{-\infty}^\infty
\left|Ve^{i\tau\sqrt{G}}\psi(h^2G)f(x)
\right|\left|g(y)\right|d\tau 
\int_{-\infty}^\infty\langle\tau/h\rangle^{1/2}\left|\widetilde 
K_h(|x-y|,\tau)\right|d\tau dxdy$$
$$\le Ch^{1/2}|t|^{-1}\int_{{\bf R}^2}\int_{{\bf R}^2}\int_{-\infty}^\infty
|\tau|^{1/2}\left\langle|x-y|/h\right\rangle^{-1/2}
\left|Ve^{i\tau\sqrt{G}}\psi(h^2G)f(x)
\right|\left|g(y)\right|d\tau dxdy$$
$$+ Ch|t|^{-1}\int_{{\bf R}^2}\int_{{\bf R}^2}\int_{-\infty}^\infty
\left|Ve^{i\tau\sqrt{G}}\psi(h^2G)f(x)
\right|\left|g(y)\right|d\tau dxdy$$
$$\le Ch^{1/2}|t|^{-1}\|f\|_{L^1}\|g\|_{L^1},$$
where $\widetilde K_h(|x-y|,t)$ denotes the kernel of the operator
$\sin\left(t\sqrt{G_0}\right)
\widetilde\psi_1(h^2G_0)$, and we have used (2.19) together with the fact that the function 
$\widetilde K_h(\sigma,t)$ satisfies (2.27). Thus, we obtain 
$$\left|\left\langle\Psi(t;h)f,g\right\rangle\right|\le Ch^{1/2}
\left\|\Psi(t;h)f\right\|_{L^\infty}\|g\|_{L^1}
+ Ch^{1/2}|t|^{-1}\|f\|_{L^1}\|g\|_{L^1},$$
which clearly implies (3.1), provided $h$ is taken small enough.
\eproof

Universit\'e de Nantes,
 D\'epartement de Math\'ematiques, UMR 6629 du CNRS,
 2, rue de la Houssini\`ere, BP 92208, 44332 Nantes Cedex 03, France

e-mail: simon.moulin@math.univ-nantes.fr


\begin{thebibliography}
\frenchspacing \baselineskip=12 pt plus 1pt minus 1pt

\bibitem{kn:B} {\sc M. Beals}, {\em Optimal $L^\infty$ decay estimates for solutions to the wave equation with a potential},
Commun. Partial Diff. Equations {\bf 19} (1994), 1319-1369.

\bibitem{kn:CCV} {\sc F. Cardoso, C. Cuevas and G. Vodev},
{\em Dispersive estimates of solutions to the wave equation with a
potential in dimensions two and three}, Serdica Math. J. {\bf 31}
(2005), 263-278.

\bibitem{kn:DP} {\sc P. D'ancona and V. Pierfelice}, {\em On the wave equation with a large rough potential},
J. Funct. Analysis {\bf 227} (2005), 30-77.

\bibitem{kn:GV} {\sc V. Georgiev and N. Visciglia}, {\em Decay estimates for the wave equation with potential},
Commun. Partial Diff. Equations {\bf 28} (2003), 1325-1369.

\bibitem{kn:G} {\sc M. Goldberg},
{\em Dispersive bounds for the three dimensional Schr\"odinger
equation with almost critical potentials}, GAFA {\bf 16} (2006),
517-536.

\bibitem{kn:JSS} {\sc J.-L. Journ\'e, A. Sofer and C. Sogge},
{\em Decay estimates for Schr\"odinger operators}, Commun. Pure
Appl. Math. {\bf 44} (1991), 573-604.

\bibitem{kn:M} {\sc S. Moulin},
{\em Low frequency dispersive estimates for the wave equation 
in higher dimensions}, submitted.

\bibitem{kn:MV} {\sc S. Moulin and G. Vodev},
{\em Low frequency dispersive estimates for the Schr\"odinger group
in higher dimensions}, Asymptot. Anal., to appear.

\bibitem{kn:S} {\sc W. Schlag},
{\em Dispersive estimates for Schr\"odinger operators in two
dimensions}, Commun. Math. Phys. {\bf 257} (2005), 87-117.

\bibitem{kn:RS} {\sc I. Rodnianski and W. Schlag},
{\em Time decay for solutions of Schr\"odinger equations with rough and time-dependent potentials}, Invent. Math. {\bf 155} (2004),
451-513.

\bibitem{kn:V1} {\sc G. Vodev},
{\em Dispersive estimates of solutions to the Schr\"odinger equation
in dimensions $n\ge 4$}, Asymptot. Anal. {\bf 49} (2006), 61-86.

\bibitem{kn:V2} {\sc G. Vodev},
{\em Dispersive estimates of solutions to the wave equation with a
potential in dimensions $n\ge 4$}, Commun. Partial Diff. Equations
{\bf 31} (2006), 1709-1733.


\end{thebibliography}
\end{document}